\documentstyle[12pt, psfig, fullpage]{amsart}
\begin{document}
\begin{abstract}
In this paper we correct an analysis of the two-player
perfect-information game Dukego given in Berlekamp, Conway, and
Guy~\cite[Chapter 19]{WW}. In particular, we characterize the board
dimensions that are fair, i.e., those for which the first player to
move has a winning strategy.
\end{abstract}
\newcommand{\scroll}[1]{\scrollmode#1\errorstopmode}
\newcommand{\comment}[1]{}

\newcommand{\createtitle}[2]{\title{#1}\author{Greg Martin}\address{Department of Mathematics\\University of Toronto\\Canada M5S 3G3}\email{gerg@@math.toronto.edu}\subjclass{#2}\maketitle}
\newcommand{\newmake}[2]{\section{#2}\label{#1sec}\noindent\input{#1}}
\newcommand{\fakemake}[1]{\section{{\tt#1.tex}}\label{#1sec}\noindent\input{#1}} %\vfil\eject}
\newcommand{\make}[1]{\label{#1sec}\noindent\input{#1}}

\newcommand{\proofmode}{\let\lebal=\label\renewcommand{\label}[1]{\lebal{##1}\quad\fbox{\smaller\smaller\tt ##1}}\hfuzz=50pt}

%%% use for AMSTeX
% \input{/opt/texmf/tex/ams/amstex}
%%%

%\comment{
\newtheorem{theorem}{Theorem}
\newtheorem{lemma}[theorem]{Lemma}
\newtheorem{corollary}{Corollary}[theorem]
\newtheorem{proposition}[theorem]{Proposition}
%} % end comment

\newenvironment{pflike}[1]{\noindent{\bf #1}}{\vskip10pt} % NO \QED
\newenvironment{proof}{\begin{pflike}{Proof:}}{\qed\end{pflike}}

\def\(#1\)#2{
 \def\helpera{\ifcase#2(\or\big(\or\Big(\or\bigg(\or\Bigg(\else(\fi}
 \def\helperb{\ifcase#2)\or\big)\or\Big)\or\bigg)\or\Bigg)\else)\fi}
 \helpera{#1}\helperb}

\newcommand{\2}[1]{\ifmmode{\cal#1}\else$\cal#1$\fi}
\newcommand{\3}[1]{\#\{#1\}}
\newcommand{\abs}[1]{\left|#1\right|} % take out
\newcommand{\floor}[1]{\lfloor#1\rfloor}
\newcommand{\bfloor}[1]{\big\lfloor#1\big\rfloor}
\newcommand{\bbfloor}[1]{\bigg\lfloor#1\bigg\rfloor}
\newcommand{\ceil}[1]{\lceil#1\rceil}
\newcommand{\bceil}[1]{\big\lceil#1\big\rceil}
\newcommand{\bbceil}[1]{\bigg\lceil#1\bigg\rceil}

\renewcommand{\mod}[1]{{\ifmmode\text{\rm\ (mod~$#1$)}
 \else\discretionary{}{}{\hbox{ }}\rm(mod~$#1$)\fi}}
\newcommand{\ep}{\varepsilon}
\renewcommand{\implies}{\Rightarrow}
\newcommand{\rmif}{\hbox{if\ }}

\newcommand{\half}{{\mathchoice{\textstyle\frac12}{1/2}{1/2}{1/2}}}
\newsymbol\dnd 232D % see The Joy of TeX (2nd ed.), Appendix G
\newcommand{\exdiv}{\mathrel{\|}}

\renewcommand{\lg}[1]{\mathop{\log_{#1}}}
\def\lgs#1^#2{\mathop{\log_{#1}^{#2}}}
\newcommand{\li}{\mathop{\rm li}\nolimits}

\newcommand{\doublespace}{%\parskip 0.2cm
  \baselineskip=24pt}
\newcommand{\spaceandahalf}{\parskip6pt\baselineskip=18pt}

\vfuzz=2pt % prohibits overfull \vbox messages due to page headers
	   % and footers

\createtitle{Restoring Fairness to Dukego}{90D46}
%\doublespace

%\renewcommand{\psfig}[1]{} %%% draft mode

\newcommand{\by}{\ifmmode\times\else$\times$\fi}

\newcommand{\D}{${\cal D}$}
\newcommand{\G}{${\cal G}$}
\newcommand{\bA}{{\bf A}}
\newcommand{\bB}{{\bf B}}
\newcommand{\bC}{{\bf C}}
\newcommand{\bD}{{\bf D}}
\newcommand{\bE}{{\bf E}}
\newcommand{\bF}{{\bf F}}
\newcommand{\bG}{{\bf G}}
\newcommand{\bH}{{\bf H}}

\newcommand{\fixfig}[2]{\vskip#1in\centerline{#2}\vskip-#1in}

\section{Introduction}\noindent

The game of Quadraphage, invented by R.~Epstein (see~\cite{Epstein}
and~\cite{Gardner}), pits two players against each other on a (generalized)
$m\by n$ chess board. The Chess player possesses a single chess piece such as a
King or a Knight, which starts the game on the center square of the board (or as
near as possible if $mn$ is even); his object is to move his piece to any square
on the edge of the board. The Go player possesses a large number of black
stones, which she can play one per turn on any empty square to prevent the chess
piece from moving there; her object is to block the chess piece so that it cannot
move at all. (Thus the phrase ``a large number'' of black stones can be
interpreted concretely as $mn-1$ stones, enough to cover every square on the
board other than the one occupied by the chess piece.) These games can also be
called Chessgo, or indeed Kinggo, Knightgo, etc.~when referring to the game
played with a specific chess piece.

Quadraphage can be played with non-conventional chess pieces as well; in fact,
if we choose the chess piece to be an ``angel'' with the ability to fly to any
square within a radius of 1000, we encounter J.~Conway's infamous
angel-vs.-devil game~\cite{AD}. In this paper we consider the case where the
chess piece is S.~Golomb's Duke, a Fairy Chess piece that is more
limited than a king, in that it moves one square per turn but only in a vertical
or horizontal direction. In this game of Dukego, we will call the Chess player
\D\ and the Go player \G.

Berlekamp, Conway, and Guy analyzed the game of Dukego in~\cite{WW}, drawing upon
strategies developed by Golomb. As they observe, moving first is never a
disadvantage in Dukego; therefore for a given board size $m\by n$, either the
first player to move has a winning strategy (in which case the $m\by n$ board is
said to be {\it fair\/}), or else \D\ has a winning strategy regardless of who
moves first, or \G\ has a winning strategy regardless of who moves first.
In~\cite{WW} it is asserted that all boards of dimensions $8\by n$ ($n\ge8$) are
fair, while \D\ has a winning strategy on boards of dimensions $7\by n$ even if
\G\ has the first move. This is not quite correct, and the purpose of this paper
is to completely characterize the fair boards for Dukego. The result of
the analysis is as follows:

\bigskip
{\narrower\noindent\it The only fair boards for Dukego are the $8\by8$ board,
the $7\by8$ board, and the $6\by n$ boards with $n\ge9$. On a board smaller
than these, \D\ can win even if \G\ has the first move, and on a board larger
than these, \G\ can win even if \D\ has the first move.

}
\bigskip
We make the convention throughout this paper that when stating the dimensions
$m\by n$ of a board, the smaller dimension is always listed first. With this
convention, the winner of a well-played game of Dukego is listed in the table
below: the entry $\ast$ denotes a fair board, on which the first player to move
can win, while the entries ``D'' and ``G'' denote boards on which the
corresponding player always has a winning strategy independent of the player to
move first.

\bigskip
\centerline{
\begin{tabular}{c||c@{\;}c@{\;}c@{\;}c@{\;}c@{\;}c@{\;}c@{\;}c@{\;}c}
% $m\by n$
& $n\le5$ &\vline& $n=6$ &\vline& $n=7$ &\vline& $n=8$ &\vline&
$n\ge9$
\\
\hline \hline
$m\le5^{\vphantom{(^)}}$ & D &\vline& D &\vline& D &\vline& D &\vline& D \\
\hline
$m=6^{\vphantom{(^)}}$ & &\vline& D &\vline& D &\vline& D &\vline& $\ast$ \\
\cline{1-1}\cline{3-10}
$m=7^{\vphantom{(^)}}$ & && &\vline& D &\vline& $\ast$ &\vline& G \\
\cline{1-1}\cline{5-10}
$m=8^{\vphantom{(^)}}$ & && && &\vline& $\ast$ &\vline& G \\
\cline{1-1}\cline{7-10}
$m\ge9^{\vphantom{(^)}}$ & && && && &\vline& G \\
\end{tabular}
}
\bigskip

As a variant of these Quadraphage games, we can allow \G\ to have both white
(wandering) and black (blocking) stones, where the white stones can be
moved from one square of the board to another once played. In this variant,
\G\ has the following options on each of her turns: place a stone of either color
on an empty square of the board, move a white stone from one square of the board
to any empty square, or pass. With a limited number of stones, it might be
the case that \G\ cannot completely immobilize the Duke, yet can play in such
a way that the Duke can never reach any of the edge squares. For instance, it
is shown in~\cite{WW} that \G\ can win (in this sense of forcing an
infinite draw) against \D\ on an 8\by8 board with only three white stones, or
with two white stones and two black stones, or with one white stone and four
black stones. In this paper we show the following:

\bigskip
{\narrower\noindent\it If \G\ has at most two white stones and at most one
black stone, then \D\ can win this variant of Dukego on a board of any size,
regardless of who has the first move. On the other hand, if \G\ has at least
three white stones, or two white stones and at least two black stones, then
the winner of this variant of Dukego is determined by the size of the board and
the player with the first move in exactly the same way as in the
standard version of Dukego (as listed in the table above).
\bigskip

}

In the analysis below, we consider the longer edges of the board to be oriented
horizontally, thus defining the north and south edges of the board, so that an
$m\by n$ board (where by convention $m\le n$) has $m$ rows and $n$ columns. Also,
if one or both of the dimensions of the board are even, we make the convention
that the Duke's starting position is the southernmost and easternmost of the
central squares of the board.

\section{How \D\ can win}\noindent

In this section we describe all the various situations (depending on the board
size, the player to move first, and the selection of stones available to \G) in
which the chess player \D\ has a winning strategy.

We begin with the simple observation that if the Duke is almost at the edge
of the board---say, one row north of the southernmost row---and \G\ has (at
most) one stone in the southernmost two rows, then \D\ can win even if it is
\G's turn. For after \G\ plays a second stone, one of the stones must be
directly south of the Duke to prevent an immediate win by \D. By symmetry,
we can assume that the second stone is to the west of the Duke, whereupon \D\
simply moves east every turn; even if \G\ continues to block the southern
edge by placing stones directly south of the Duke on each turn, \D\ will
eventually win by reaching the east edge. If \D\ is in this situation, we say
that \D\ has an Imminent Win on the south edge of the board (and similarly
for the other edges).

It is now easy to see that \D\ can always win on a $5\by n$ board even if
\G\ has the first move. The Duke will be able to move either north or south
(towards one the long edges) on his first turn, since \G\ cannot block both
of these squares on her first turn; and then \D\ will have an Imminent Win on
the north or south edge of the board, correspondingly. Of course,
this implies that \D\ can win on any $3\by n$ or $4\by n$ board as
well, regardless of who moves first (Dukego on $1\by n$ and $2\by n$
boards being less interesting still). Similarly, \D\ can win on a $6\by n$
board if he has the first move, since he can move immediately into an
Imminent Win situation along the south edge of the board (recalling our
convention that the Duke starts in the southernmost and/or easternmost
central square).

We can also see now that \D\ can win on any size board if \G\ is armed with
only two white (wandering) stones. \D\ selects his favorite of the four
compass directions, for instance south, and pretends that the row directly
adjacent to the Duke in that direction is the edge of the board. By adopting an
Imminent Win strategy for this fantasy edge row, \D\ will be able either to
reach the east or west edge of the board for a true win, or else to move one row
to the south for a fantasy win. But of course, repeating this Fantastic Imminent
Win strategy will get the Duke closer and closer to the south edge of the board,
until his last fantasy win is indeed a win in reality. Similarly, \D\ can
always win if \G\ adds a single black (blocking) stone to her two white stones.
\D\ plays as above until \G\ is forced to play her black stone (if she never
does, then we have just seen that \D\ will win); once the black stone is played,
\D\ rotates the board so that the black stone is farther north than the Duke,
and then uses this Fantastic Imminent Win strategy towards the south edge.

\begin{figure}[b]
\vskip4.7in
\centerline{\psfig{figure=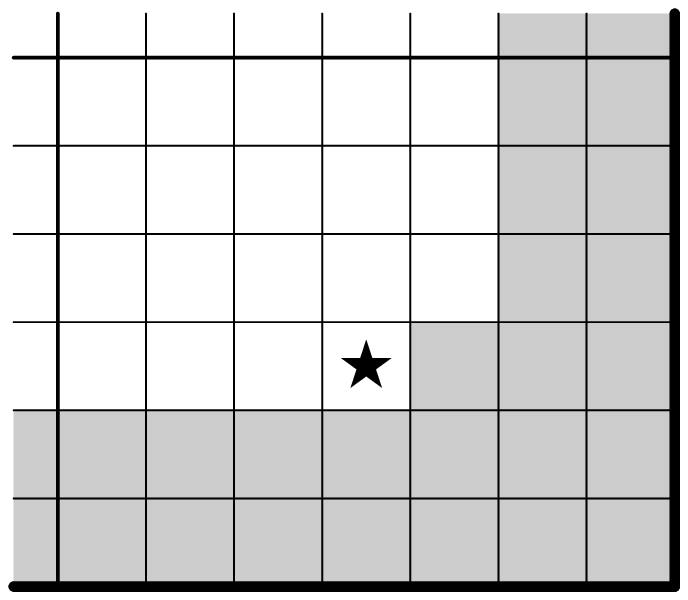}\hfil\hfil\raise15mm\psfig{figure=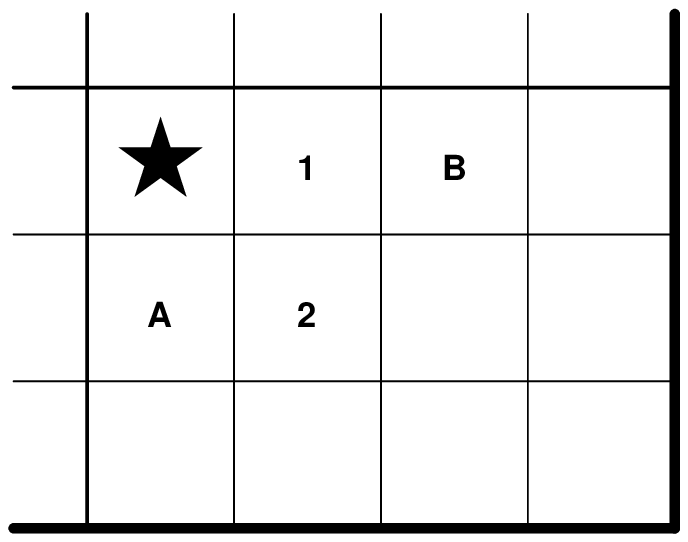}}
\vskip-4.7in
\caption{The Corner Win: If \G\ has no stones in the shaded region
(left), then \D\ can win from position $\star$ (right)}
\label{corner}
\end{figure}
We now describe a slightly more complicated situation in which \D\ has
a winning strategy. Suppose that the Duke is two rows north of the
southernmost row and three rows west of the easternmost row, and that
there are no stones in the two southernmost rows or in the two
easternmost rows of the board, nor is there a stone directly east of
the Duke (see the left diagram of Figure~\ref{corner}).  Then we claim
that \D\ can win even if it is \G's turn to move. For (referring to
the second diagram of Figure~\ref{corner}) \G\ must put a stone in
square \bA\ to block \D\ from moving into an Imminent Win along the
south edge. \D\ moves east to square {\bf1}, whereupon \G\ must
put a stone in square \bB\ to block an Imminent Win along the east
edge. But this is to no avail, as \D\ then moves to square {\bf2} to
earn an Imminent Win along the south edge anyway. If \D\ is in this situation
described in Figure~\ref{corner}, we say that \D\ has a Corner Win in
the southeast corner of the board (and similarly for the other corners).

Notice that it is necessary for the shaded area to be empty for the
Corner Win to be in force. If \G\ has a stone in the
second-southernmost row somewhere to the west of the Duke, then she
can place a stone on square {\bf2} to safely guard the south edge of
the board. Similarly, if \G\ has a stone in the second-easternmost row
somewhere to the north of the Duke, then she can defend both edges of
the board by playing a stone at square A on her first turn and one at
square {\bf2} on her second turn.

We can now see that \D\ can win on a 6\by8 board regardless of who has
the first move, since the Duke starts the game in a Corner Win
position. This implies that \D\ can win on 6\by6 and 6\by7 boards
regardless of who has the first move. Also, we can argue that \D\ can
win on a 7\by7 board even if \G\ has the first move. By symmetry, we
can assume that \G's first stone is placed to the northwest of the
Duke or else directly north of the Duke, whereupon \D\ can move to the
south and execute a Corner Win in the southeast corner of the
board. Similarly, \D\ can win on an 8\by8 board (and thus on a 7\by8
board as well) if he has the first move, since he can again gain a
Corner Win situation by moving south on his first turn.

\section{How \G\ can win}

We have now shown all we claimed about \D's ability to win; it's time
to give \G\ her turn.
\begin{figure}[b]
\vskip4.9in
\centerline{\psfig{figure=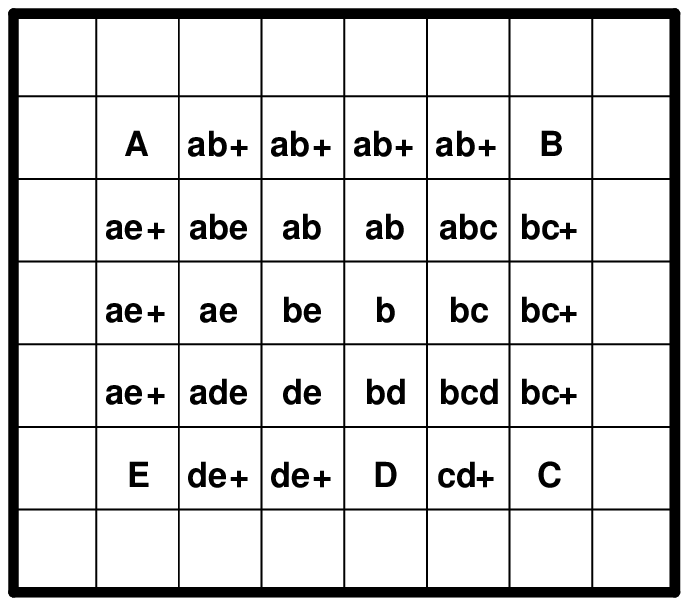}\hfil\hfil\raise22mm\psfig{figure=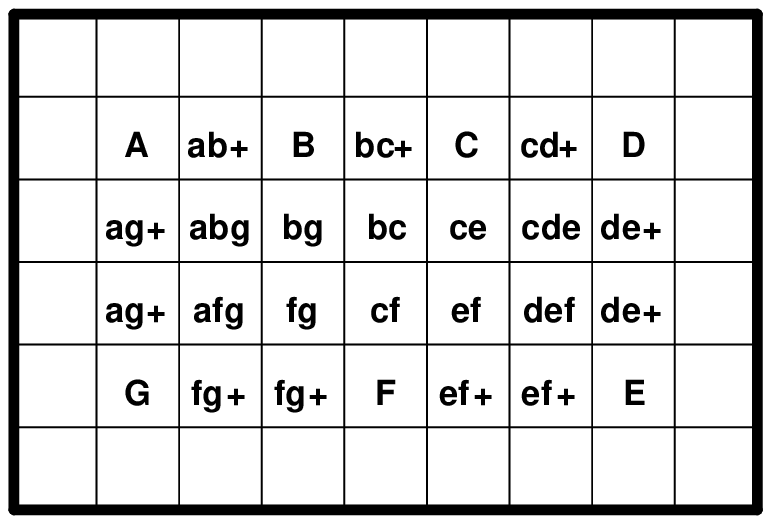}}
\vskip-4.9in
\caption{\G's strategies on a 7\by8 board (left) and on a 6\by9 board
(right) with three white stones}
\label{strategies}
\end{figure}
We start by exhibiting strategies for \G\ to win on a 7\by8 board and on a 6\by9
board with the first move. To begin, we assume that \G\ possesses only three
white (wandering) stones; since having extra stones on the board is never
disadvantageous to \G, the strategy will also show that \G\ can win with a large
number of black stones and no white stones. 

The squares labeled with capital letters in Figure~\ref{strategies}
are the strategic squares for \G's strategies on these boards. The key
to reading the strategies from Figure~\ref{strategies} is as follows:
whenever the Duke moves on a square labeled with some lowercase
letters, \G\ must choose her move to ensure that the squares with the
corresponding capital letters are all covered. Other squares may be
covered as well, as this is never a liability for \G. If the Duke is
on a square with a $+$ sign as well as some lowercase letters, \G\
must position a tactical stone on the edge square adjacent to the Duke
(blocking an immediate win) as well as having strategic stones on the
corresponding uppercase letters.  All that is required to check that
these are indeed winning strategies is to verify that every square
marked with lowercase letters (counting $+$ as a lowercase letter for
this purpose) contains all of the letters, save at most one,
of each of its neighbors, so that \G\ can correctly change
configurations by moving at most one white stone.

In the case of the 7\by8 board, the Duke begins on the square marked
``b'' in the left-hand diagram of Figure~\ref{strategies} (recalling
our convention about the precise beginning square for \D\ on boards
with one or both dimensions even), so \G's first move will be to place
a stone on the square marked \bB. Notice how \G\ counters the instant
threat of a Corner Win by \D\ in the southeast corner, by playing her
first stone on the second-easternmost column of the board at \bB, a
move which also begins the defense of the north and east edges of the
board against direct charges by the Duke.

In the case of the 6\by9 board, we need to make another convention
about the beginning of the game. The Duke begins on the more southern
of the two central squares (the square marked ``cf'' in the right-hand
diagram of Figure~\ref{strategies}), and we stipulate that \G's first
stone be played at square \bF. We may also assume that \D's first move
is not to the north, for in this case \G\ may rotate the board 180
degrees and pretend that \D\ is back in the starting position, moving
her white stone from the old square \bF\ to the new square
\bF. Eventually, \D\ will move east or west, to a square marked either
``ef'' or ``fg'', and at this point \G\ begins to consult the
right-hand diagram for her strategy, playing her second stone at
square \bE\ or \bG, respectively.

\G's strategy on the 7\by8 board can be converted into a strategy
using two white stones and two black stones without too much
difficulty. As before, \G\ begins by placing a white stone on the
square \bB. \G's goal is to establish her two black stones either on
squares \bA\ and \bC\ (or on \bB\ and \bE), and then use her white stones
both on the strategic squares \bB\ and \bE\ (or \bA\ and \bC,
respectively) and as tactical stones blocking immediate wins. The
strategic square \bD\ is only used to keep \D\ in check until the two
black stones can be established. The conversion is straightforward and we
omit the details.

\begin{figure}[btf]
\vskip4.1in
\centerline{\psfig{figure=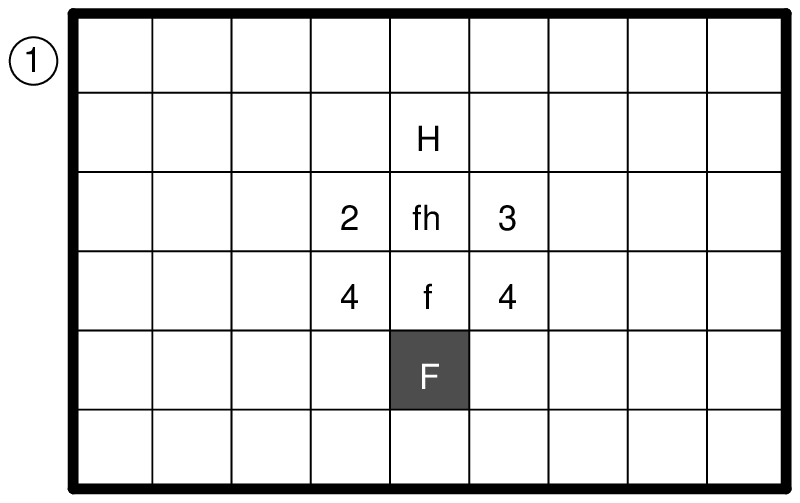}}
\centerline{\psfig{figure=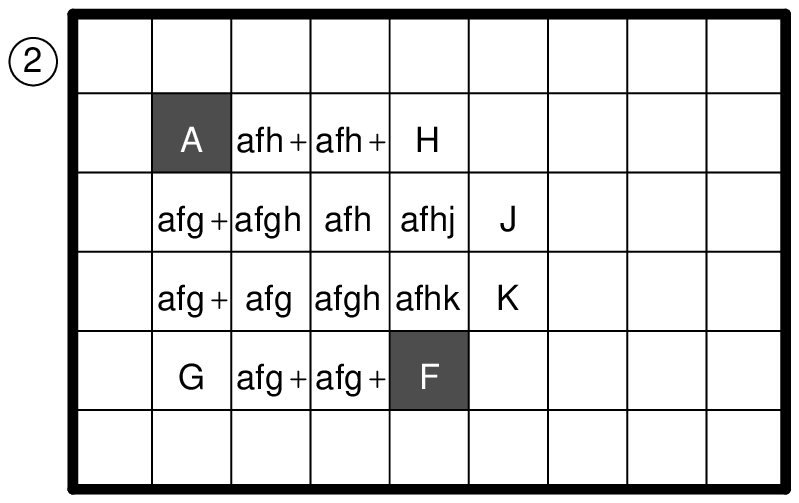}\hfil\hfil\psfig{figure=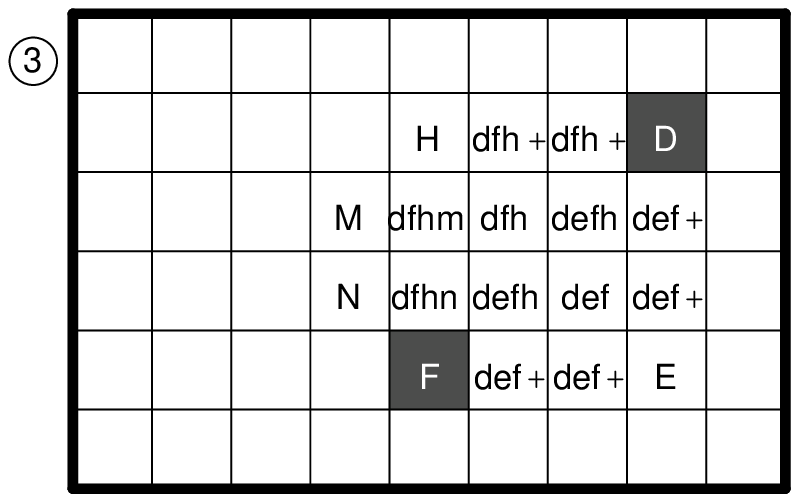}}
\centerline{\psfig{figure=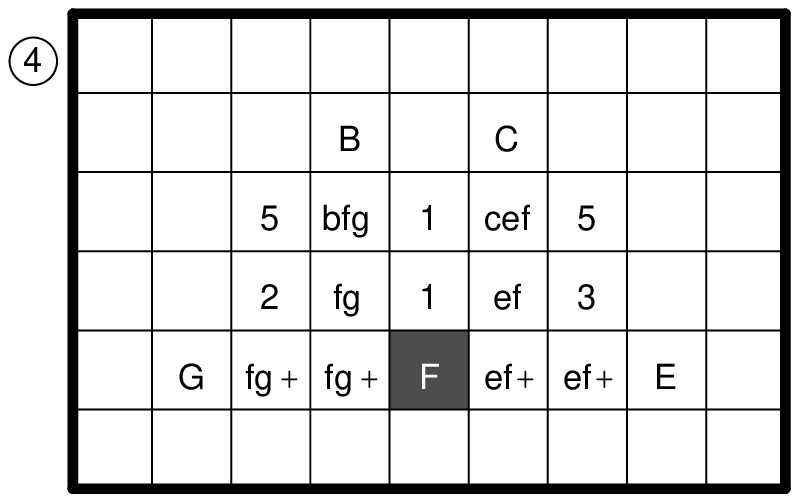}\hfil\hfil\psfig{figure=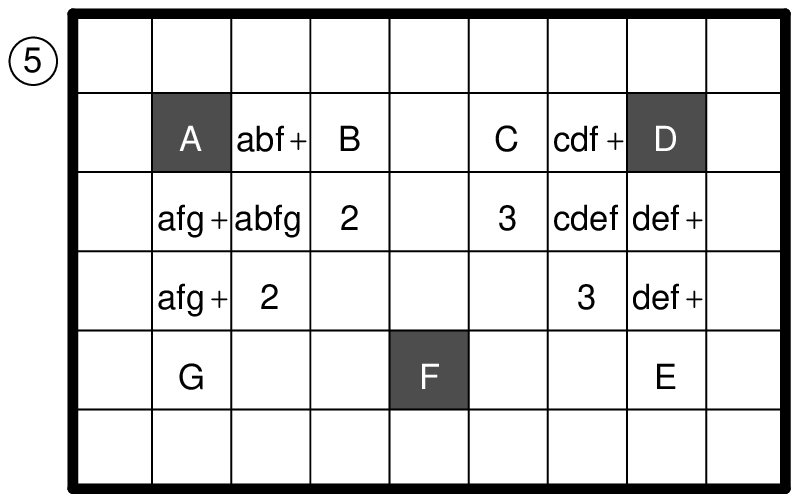}}
\vskip-4.1in
\caption{\G's strategy on a 6\by9 board with two white stones and two
black stones}
\label{strategy69b}
\end{figure}
Less straightforward, however, is showing that \G\ also has a strategy for
winning on the 6\by9 board with the first move, if she has two white stones and
two black stones. We include in Figure~\ref{strategy69b} a full strategy for \G\
in this situation. \G's goal is to coerce \D\ into committing to one of the two
sides of the board, corresponding to diagram 2 or diagram 3 in
Figure~\ref{strategy69b}.  Then she will be able to establish her black stones
on squares \bA\ and \bF\ (or \bD\ and \bF), and use her white stones both as
strategic stones on squares \bG\ and \bH\ (or \bE\ and \bH, respectively) and as
tactical stones blocking immediate wins for \D.

\G\ begins by reading the topmost diagram (labeled 1). Since the Duke starts on
the square marked ``f'', \G\ places a stone on the square \bF; since square \bF\
is shaded black in the diagram, this stone must be a black stone. When the Duke
moves to a square marked with a number, \G\ immediately switches to the
corresponding diagram in Figure~\ref{strategy69b} and moves according to the
Duke's current position. For example, suppose that the Duke's first move is to
the east, onto the square marked {\bf4}; in this case \G\ switches to diagram 4,
where the Duke's square is marked ``ef'', indicating that \G\ must add a stone
to square \bE\ (a white stone, since square \bE\ is not shaded black) in
addition to her existing black stone on square \bF. We remark, to
ameliorate one potential source of confusion, that diagram 5 is really
a combination of two smaller diagrams, one for each side of the board;
in particular, there will never be a need to have black stones
simultaneously at squares \bA, \bD, and \bF.

Of course, any opening move for \G\ other than placing a stone on square \bF\
would lead to a quick Imminent Win for \D\ along the south edge. It turns out,
though, that even if \G's first move is to play a {\it white\/} stone at \bF, a
winning strategy exists for \D\ (assuming still that \G\ has exactly two
stones of each color at her disposal). A demonstration of this would be somewhat
laborious, and so we leave the details for the reader's playtime.

To conclude this section, we are now in a position to argue that \G, armed with
a large number of black stones or with three white stones or with two stones of
each color, can win on a 7\by9 board (and thus on any larger board) even if \D\
has the first move. By symmetry we may assume that \D's first move is either to
the south or to the east. If \D\ moves east on his first move, then \G\ ignores
the westernmost column of the board and adopts the above 7\by8 strategy on the
remainder of the board; alternatively, if \D\ moves south on his first turn,
then \G\ ignores the northernmost row of the board and adopts the above 6\by9
strategy on the rest of the board.

\section{Afterthoughts}

It is not quite true that we have left no stone unturned (pun unintended) in
our analysis of Dukego. For instance, it is unclear exactly how many black
stones \G\ needs to win the original version of Dukego (no white stones) on the
various board sizes; determining these numbers would most likely involve a fair
amount of computation to cover \D's initial strategies.

Somewhat more tractable, however, would be to determine how many black stones
\G\ needs to win when she also possesses a single white stone. The strategy
given in~\cite{WW} for \G\ on an 8\by8 board works perfectly well when \G\ has
one white stone (to be used tactically) and four black stones (to be placed
strategically), and this is the best that \G\ could hope for. On the other hand,
\G's strategies on the 7\by8 and 6\by9 boards as described above require five
and seven black stones, respectively, to go along with the single white stone
(and it requires some care to show that seven black stones suffice for the
6\by9 board, since we need to account for \D\ moving north on his first move).

We believe that \G\ {\it cannot\/} win on either a 7\by8 board or a 6\by9 board
with a single white stone and only four black stones. If this is the case, then
the least number of black stones that \G\ can win with, when supplemented by a
white stone, would be five on a 7\by8 board; but we don't know whether the
analogous number on a 6\by9 board is five, six, or seven.

\bibliography{dukego}
\bibliographystyle{../amsplain}
\end{document}